\documentclass[11pt]{article}

\usepackage{epsfig }
\usepackage{latexsym}
\usepackage{amsfonts,amsmath,amssymb,textcomp}
\allowdisplaybreaks

\newtheorem{thm}{Theorem}[section]
\newtheorem{lem}[thm]{Lemma}
\newtheorem{cor}[thm]{Corollary}
\newtheorem{prop}[thm]{Proposition}
\newtheorem{rem}[thm]{Remark}
\newtheorem{deff}[thm]{Definition}
\newtheorem{conj}[thm]{Conjecture}
\newtheorem{key}[thm]{Keywords}
\newtheorem{prob}[thm]{Problem}
\newenvironment{probb}{ \begin{prob} \rm}{ \end{prob} }
\newcommand{\bprobb}{\begin{probb}}
\newcommand{\eprobb}{\end{probb}}
\newcommand{\bth}{\begin{thm}}
\newcommand{\ethGL}{\end{thm}}
\newcommand{\bconj}{\begin{conj}}
\newcommand{\econj}{\end{conj}}
\newcommand{\bkey}{\begin{key}}
\newcommand{\ekey}{\end{key}}
\newcommand{\bl}{\begin{lem}}
\newcommand{\el}{\end{lem}}
\newcommand{\bdeff}{\begin{deff}}
\newcommand{\edeff}{\end{deff}}
\newcommand{\bcor}{\begin{cor}}
\newcommand{\ecor}{\end{cor}}
\newcommand{\bprop}{\begin{prop}}
\newcommand{\eprop}{\end{prop}}
\newcommand{\brem}{\begin{rem}}
\newcommand{\erem}{\end{rem}}
\newcommand{\beq}{\begin{equation}}
\newcommand{\eeq}{\end{equation}}
\newcommand{\beqn}{\begin{eqnarray}}
\newcommand{\eeqn}{\end{eqnarray}}
\newcommand{\beqns}{\begin{eqnarray*}}
\newcommand{\eeqns}{\end{eqnarray*}}
\newcommand{\ba}{\begin{array}}
\newcommand{\ea}{\end{array}}
\newcommand{\bit}{\begin{itemize}}
\newcommand{\eit}{\end{itemize}}
\newcommand{\ben}{\begin{enumerate}}
\newcommand{\een}{\end{enumerate}}

\newcommand{\BO}{\mathcal{O}}
\newcommand{\BN}{\mathcal{N}}

\newcommand{\babs}{\begin{abstract}}
\newcommand{\eabs}{\end{abstract}}
\newcommand{\bal}{\begin{align}}
\newcommand{\bals}{\begin{align*}}
\newcommand{\bs}{\begin{skip}}
\newcommand{\eal}{\end{align}}
\newcommand{\eals}{\end{align*}}
\newcommand{\es}{\end{skip}}

\newcommand{\ra}{\rightarrow}

\newcommand{\B}{\quad}

\newcommand{\lp}{\left (}
\newcommand{\rp}{\right )}
\newcommand{\lb}{\left [}
\newcommand{\rb}{\right ]}

\newcommand{\de}{\delta}
\newcommand{\De}{\Delta}

\newcommand{\al}{\alpha}
\newcommand{\be}{\beta}

\newcommand{\Fi}{\varphi}

\newcommand{\II}{\infty}

\newcommand{\sig}{\sigma}

\newcommand{\gam}{\gamma}

\def\1{{\ifmmode 1\mskip-1.5\thinmuskip\mathrm{l}%
        \else\textrm{1\hskip -.23em l}\fi}}

\newcommand{\E}{\mbox{$\mathbb E$}} 
\newcommand{\V}{\mbox{$\mathbb V$}} 
\def\P{{\mathbb {P}}}

\title{\bf The perimeter of uniform and geometric words: a probabilistic analysis }
\author{Guy Louchard\thanks{Universit\'e Libre de Bruxelles,
D\'epartement d'Informatique, CP 212, Boulevard du Triomphe, B-1050
Bruxelles, Belgium, email: louchard@ulb.ac.be}
}
\date{\today}
\begin{document}
\maketitle

\babs
Let a word be a sequence of $n$ i.i.d.  integer random variables. The  perimeter  $P$ of the word is the number of edges of the word, seen as a polyomino. In this paper, we present 
 a probabilistic approach  to the computation of the moments of $P$. This is applied to uniform and geometric random variables. We also show that, asymptotically, 
the distribution of $P$ is Gaussian and, 
 seen as a stochastic process, the perimeter converges in distribution to a Brownian motion.
\eabs


\section{Introduction}\label{S1} 

\begin{figure}[htbp]
	\centering
		\includegraphics[width=0.6\textwidth,angle=270]{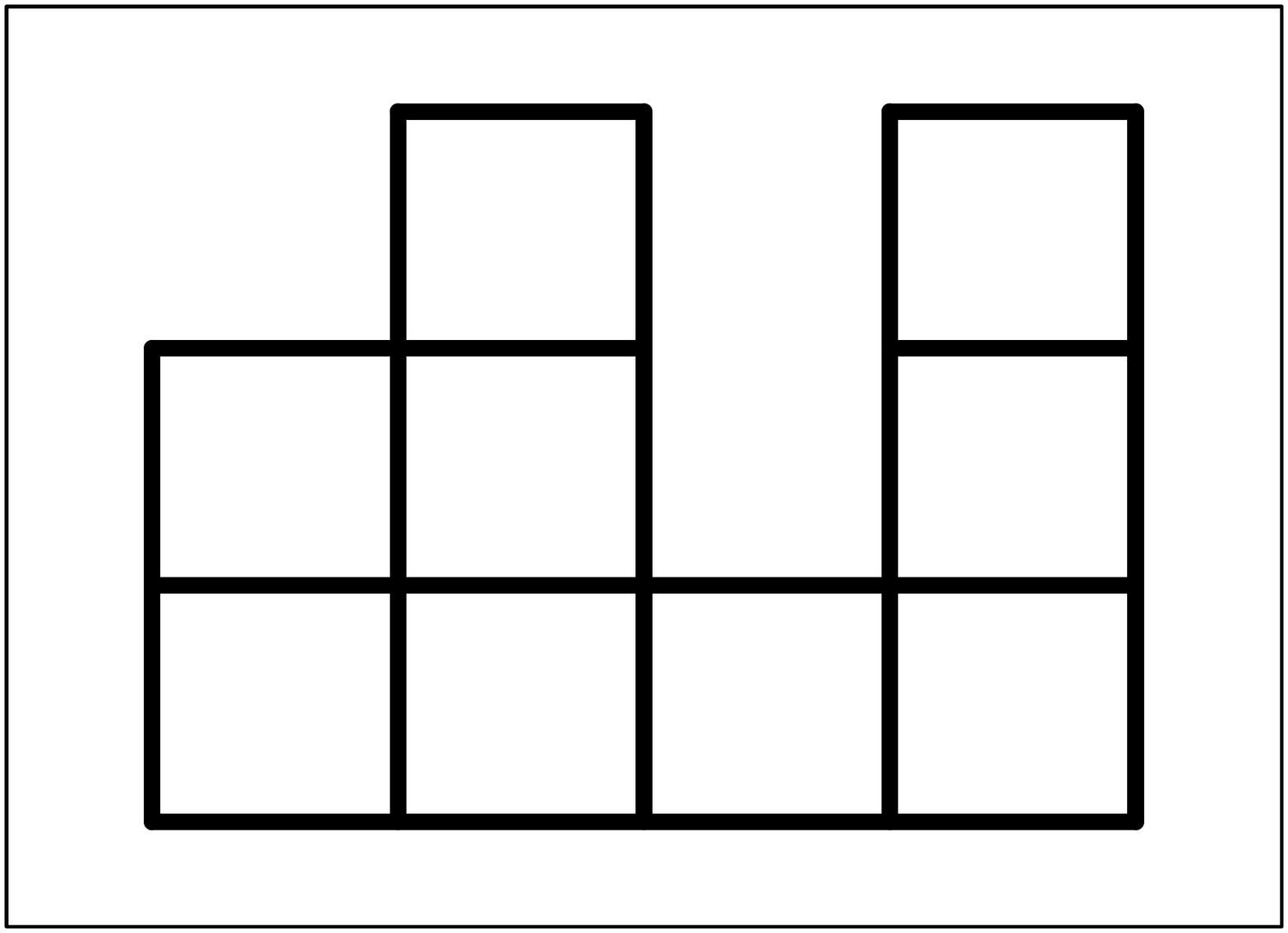}
	\caption{The polyomino based on the word ${2,3,1,3},n=4,P=18$}
	\label{F1}
\end{figure}

Our attention was recently attracted by a paper by Blecher et al. \cite{BBKM17} on the perimeter of words:
a word is a sequence of $n$ i.i.d.  integer random variables (RV) $\{x_0,x_1,\ldots,x_m\},m:=n-1$. In \cite{BBKM17}, the RV are distributed uniformly on $[1,k]$. These RV are also used in this paper. The  perimeter  $P_n$ of the word is the number of edges of the word, seen as a polyomino.
A typical polyomino, based on the word ${2,3,1,3},n=4,P=18$ is given in Fig.\ref{F1}.

 The mean $M_{P,n}:=\E(P_n)$ and variance $\V(P_n)$ of $P_n$ are given in \cite{BBKM17}, with $M=\frac{(k - 1) (k + 1)}{ 3 k} $ by the following theorem:
\bth
In the uniform $[1,k]$ case, $M_{P,n}$ and $\V(P_n)$ are given in \cite{BBKM17} by 
\bal
M_{P,n}&=(n-1)M+2n+(k+1)=\frac{(3 k + 2 k^2  + 1) + (k^2  + 6 k - 1) n}{3k},                      \label{E1}\\
\V(P_n)&=\frac{(-5 k^2  + 4 k^4  + 1) + (-3 + 3 k^4 ) n}{45k^2}.                                   \label{E2}
\end{align} 
\ethGL
Some years ago, we had been interested in some uniformly distributed words: see \cite{DL05}. Moreover, we had analyzed some polyominoes, for instance in \cite{GL96} and \cite{GL299}, where,  in Particular, we had derived some limiting Brownian motion (BM) 
Processes for trajectories. Some recent Papers on polyomino's perimeter are, for instance, \cite{BR03}, \cite{FE07}.

Another classical distribution is the classical geometric$(p)$ one, with distribution $pq^{i-1},i\geq 1,p\in(0,1),q:=1-p$. 
In several papers (some of them with H. Prodinger) we had analyzed related word parameters from a probabilistic point of view. Our last  papers on this topic being \cite{LP06}, \cite{LOPR05}. We again derived some limiting  BM processes, for instance in 
\cite{GL00}, \cite{PL01}. For other recent papers on geometric words, see \cite{ABCK17}, \cite{AK14}. 

In the present paper, our motivation is to present a novel  approach to the words perimeter problem:
\bit
\item a probabilistic approach easily leads to the moments of $P_n$,
\item the distribution of $P_n$ is asymptotically shown to be Gaussian,
\item seen as a stochastic process, the perimeter converges in distribution to a BM,
\item our technique is applied to the geometric$(p)$ case.
\eit

\section{The mean and variance of $P_n$ in the uniform $[1,k]$ case}\label{S2} 
In this section, we present a probabilistic approach to the mean and variance of the full perimeter $P_n$.

Set $Q_m:=\sum_1^m y_i, y_i:=|x_i-x_{i-1}|$. Clearly, $P_n=Q_m+x_0+x_m+2n$. For further use, we define the vertical perimeter $R_n:=Q_m+x_0+x_m$. We see that the $y_i$ are identically distributed, $y_i$ is correlated with $y_{i+1}$, but \emph{independent} of $y_k,k\geq i+2$. 

The following notations and relations  will be used throughout the paper: 
 \bals
 m&:=n-1,\\
 \overline{z}&:=z-\E(z),\mbox{ for any RV},\\
 M&:=\E(y_i),\\
 M_{Q,m}&:=\E (Q_m)=mM,\\
 M_{R,n}&:=\E(R_n)=M_{Q,m}+2\E(x_0),\\
 M_{P,n}&:=\E(P_n)=M_{R,n}+2n,\\
 T_{\al,\be, \gam,\de}&:=\E\lp x_0^\al \cdot y_1^\be\cdot y_2^\gam\cdot y_3^\de \rp,\\
 \overline{T}_{0,\be, \gam,\de}&:=\E\lp (y_1-M)^\be\cdot (y_2-M)^\gam\cdot (y_3-M)^\de \rp,
\end{align*} 
and the same definitions for $T_{\al},T_{\al,\be},T_{\al,\be, \gam}$. When an exponent is null, it means the absence of the relevant variable.  For instance, $M \equiv T_{0,1},T_1=\E(x_0)=\frac{k+1}{2}$ for the uniform case. Explicitly, we have
\[\left. T_{\al,\be, \gam,\de}=\lb \sum_{i=1}^{k}i^\al \sum_{j=1}^{k}|j-i|^\be \sum_{\ell=1}^{k}|\ell-j|^\gam \sum_{r=1}^{k}|r-\ell|^\de\rb\right/k^4.\]
Let us first compute the distribution of $y_i$: $f(u):=\P(y_i=u), u\in [0,k-1]$. Consider first the case  $u>0$. If $x_1>x_0, u=x_1-x_0,x_1=u+x_0$. But $1\leq x_1\leq k$, hence $1\leq x_0\leq k-u$. So we first have 
\[S_1:=\frac1k\sum_1^{k-u}\P(x_0=i)=\frac{k-u}{k^2}.\]

\begin{figure}[htbp]
	\centering
		\includegraphics[width=0.8\textwidth,angle=0]{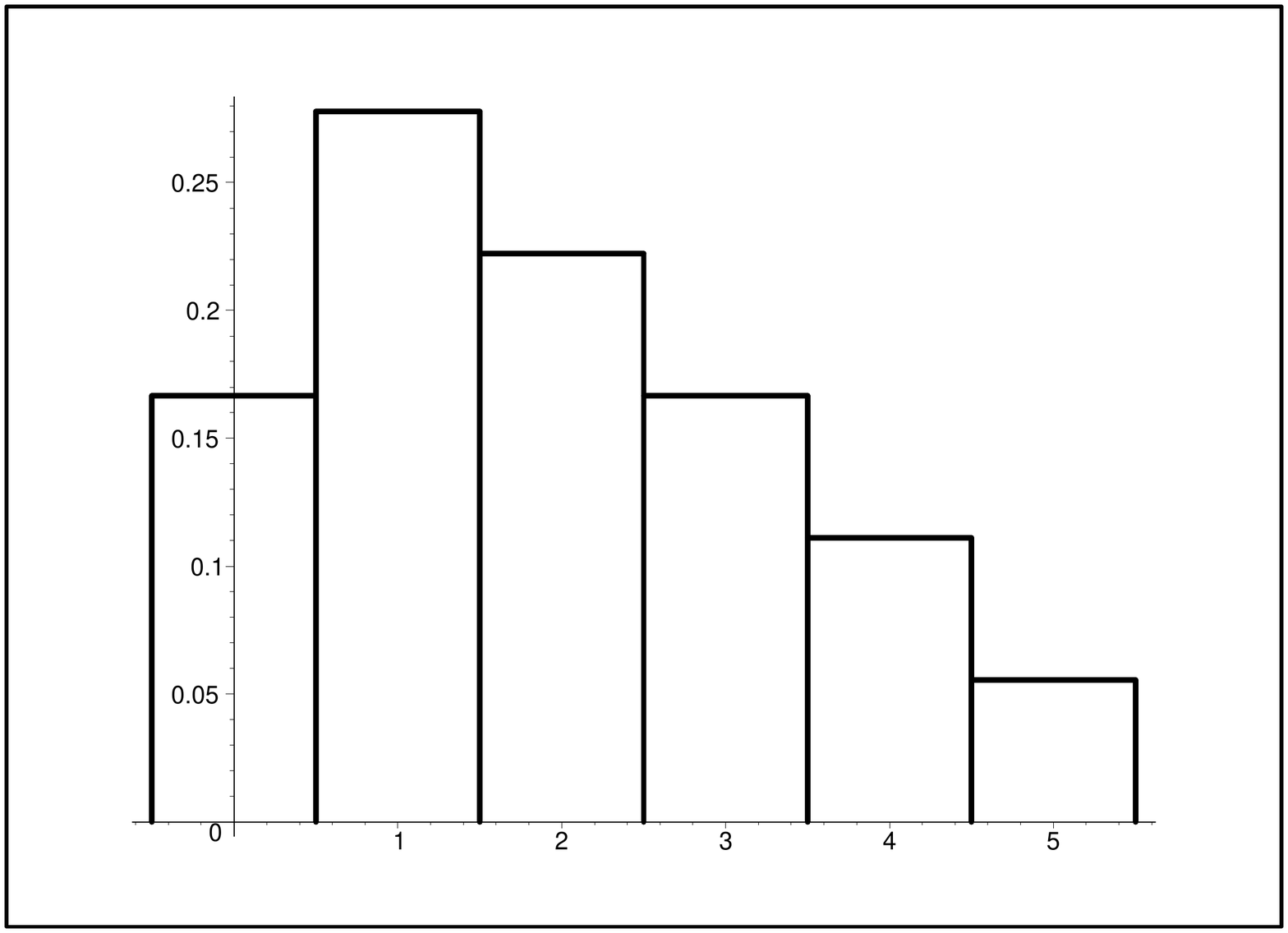}
	\caption{$f(u),k=6$ }
	\label{F2}
\end{figure}

Next, if  $x_1<x_0, u=x_0-x_1,x_1=x_0-u$, But $1\leq x_1\leq k$, hence $1+u\leq x_0\leq k$.  So 
 
\[S_2:=\frac1k\sum_{1+u}^{k}\P(x_0=i)=\frac{k-u}{k^2}.\]  
Finally,
\[f(u)=S_1+S_2=\frac{2(k-u)}{k^2},u>0. \] 
In the case $u=0$, we simply have $f(0)=\frac{1}{k^2}   \sum_1^k 1  =\frac1k$. A plot of $f(u),k=6$, is given in Fig.\ref{F2}.

Now we are ready to compute $M$. This is given either by  
\[M:=T_{0,1}:=\left. \sum_{i=1}^{k}\lp  \sum_{j=i}^{k}(j-i)+\sum_{j=1}^{i-1}(i-j)\rp \right/k^2=\frac{(k - 1) (k + 1)}{ 3 k} ,\]                   
	or $\sum_0^{k-1} f(u)u$, which of course leads to the same result.

 Hence 
  \bals    
 M_{R,n}&=(n-1)M+(k+1)=\frac{(3 k + 2 k^2  + 1) + (k^2  - 1) n}{3k},\\
 M_{P,n}&=(n-1)M+2n+(k+1)=\frac{(3 k + 2 k^2  + 1) + (k^2  + 6 k - 1) n}{3k},
	  \end{align*}                                   
which fits with(\ref{E1}).

Some useful expressions will be used in this section. We collect them here.		
\[T_1=\E(x_0)=\sum_{i=1}^{k}i/k=\frac{k+1}{2},\]	

 \[T_{2}=\E(x_0^2)=\sum_{i=1}^{k}i^2/k=\frac{   (k + 1) (1 + 2 k)}{6} ,\]

\[ T_{1,1}=\left. \sum_{i=1}^{k}i\lp \sum_{j=i}^{k}(j-i)+\sum_{j=1}^{i-1}(i-j)\rp \right/k^2=\frac{ (k - 1) (k + 1)^2}{6k},\]

\[ T_{0,2}=\left. \sum_{i=1}^{k}\sum_{j=1}^{k}(j-i)^2\right/k^2=\frac{ (k - 1) (k + 1)}{ 6} ,\mbox{ this is also given by }
\sum_{u=1}^{k-1}f(u)u^2,\]

 \[\overline{T}_{0,2}=\left. \sum_{i=1}^{k}\lp \sum_{j=i}^{k}(j-i-M)^2+\sum_{j=1}^{i-1}(i-j-M)^2\rp \right/k^2
=\frac{(k - 1) (k + 1) (k^2  + 2)}{18k^2},\]

 \bals                                
T_{0,1,1}&=\sum_{i=1}^{k}\lb \sum_{j=i}^{k}(j-i)\lp \sum_{\ell=j}^{k}(\ell-j)+\sum_{\ell=1}^{j-1}(j-\ell)\rp
 +\left. \sum_{j=1}^{i-1}(i-j)\lp \sum_{\ell=j}^{k}(\ell-j)+\sum_{\ell=1}^{j-1}(j-\ell)\rp \rb \right/k^3\\
&=\frac{(k - 1) (k + 1) (7 k^2  - 8)}{60 k^2},
 \end{align*}

 \bals                                
\overline{T}_{0,1,1}&=\sum_{i=1}^{k}\lb \sum_{j=i}^{k}(j-i-M)\lp \sum_{\ell=j}^{k}(\ell-j-M)+\sum_{\ell=1}^{j-1}(j-\ell-M)\rp \right.\\
 &\left. \left.  +\sum_{j=1}^{i-1}(i-j-M)\lp \sum_{\ell=j}^{k}(\ell-j-M)+\sum_{\ell=1}^{j-1}(j-\ell-M)\rp \rb \right/k^3\\
&=\frac{(k - 1) (k - 2) (k + 2) (k + 1)}{180 k^2}.
 \end{align*}

 These expressions are the only necessary ones in order to compute $\V(P_n)$.

	Now we turn to the computation of  variance  $\V(P_n)$. Of course, only $R_n$ has to be used here.	The dominant term of 	 $\V(R_n)$	is immediate: this is  given by
	\bals												
  &n[(T_{0,2}-M^2)+2(T_{0,1,1}-M^2)]=nV^*,\\
	V^*&=\frac{ (k - 1) (k + 1) (k^2  + 1)}{15k^2}.
	\end{align*}   
	Indeed,
	\[\overline{R}_n=\overline{x}_0+\overline{x}_m+\overline{Q}_m,\]
	and the effect of $\overline{x}_0$ on the variance is just $\overline{T}_2+2\overline{T}_{1,1}=\BO(1)$. Similarly for the
	contribution of $x_m$. Also the contribution of the couples $\overline{y}_i\overline{y}_{i+1}$ is given by
	$2(m-1)\overline{T}_{0,1,1}=2m\overline{T}_{0,1,1}+\BO(1)$ and all other contributions are null by independence. 
		 $V^*$  is of course also given by $\overline{T}_{0,2} +2\overline{T}_{0,1,1}$.

 To compute  $\V(P_n)$, we must collect all necessary terms. We symbolically expand
\[(x_0+x_m+y_1+y_i+y_{i+1}+y_{i+2}+y_m)^2.\]   
We collect the relevant contributions, with their  weights (we just have to count the corresponding tuples, and, as explained above,
all other tuples do not contribute to the variance) :
\bals
y_i^2&\ra m T_{0,2},\\
y_i y_{i+1}&\ra 2(m-1)T_{0,1,1},\\
x_0^2,x_m^2&\ra 2T_2,\\
x_0 y_1, x_m y_m &\ra 4T_{1,1},\\
y_i y_j,i+2\leq j\leq m,1\leq i\leq m-2,&\ra (m-1)(m-2)M^2, \mbox{ independent RV}\\
x_0x_m&\ra 2T_1^2,\\
x_0y_i,i>1,x_m y_j,j < m&\ra 4T_1 M(m-1), \mbox{ independent RV}.
 \end{align*}                             
 This gives                              

\bals
 \V(R_n)&=(n-1)T_{0,2}+2T_{2}+(n-2)2T_{0,1,1}+4T_{1,1}+(n-2)(n-3)M^2+2T_1^2+4T_1M(n-2)-M_{R,n}^2\\
&=\frac{(-5 k^2  + 4 k^4  + 1) + (-3 + 3 k^4 ) n}{45k^2},
\end{align*} 
which fits with (\ref{E2}).

\section{The third centered moment $\mu_3(P_n)$ of $P_n$ in the uniform $[1,k]$ case}\label{S3} 
In this section, we apply our probabilistic technique to the third centered moment computation.
We will only compute the $n-$dominant term of $\mu_3(P_n)$, the complete analysis  goes as in Sec. \ref{S2}, only with elementary but tedious algebra, we omit the details.

 The necessary expressions are as follows (for the sake of completeness, we also provide the centered moments):

 \[\overline{T}_{0,3}=\left. \sum_{i=1}^{k}\lp \sum_{j=i}^{k}(j-i-M)^3+\sum_{j=1}^{i-1}(i-j-M)^3\rp \right/k^2
=\frac{ (k - 1) (k - 2) (k + 2) (k + 1) (2 k^2  - 5)}{270k^3},\]

 \[T_{0,3}=\left. \sum_{u=1}^{k-1}f(u)u^3=\sum_{i=1}^{k}\lp \sum_{j=i}^{k}(j-i)^3+\sum_{1}^{i-1}(i-j)^3\rp \right/k^2=\frac{    (k - 1) (k + 1) (3 k^2  - 2)}{30 k},\]

\bals
\overline{T}_{0,1,1,1}&=\sum_{i=1}^{k}\lb \rule{0mm} {10mm}\right.\sum_{j=i}^{k}(j-i-M)\lb \sum_{\ell=j}^{k}(\ell-j-M)\lp \sum_{r=\ell}^{k}(r-\ell-M)
+\sum_{r=1}^{\ell-1}(\ell-r-M)\rp \right. \\
&\left.  +\sum_{\ell=1}^{j-1}(j-\ell-M)\lp \sum_{r=\ell}^{k}(r-\ell-M)+\sum_{r=1}^{\ell-1}(\ell-r-M)\rp \rb\\
& +\sum_{j=1}^{i-1}(i-j-M)\lb \sum_{\ell=j}^{k}(\ell-j-M)\lp \sum_{r=\ell}^{k}(r-\ell-M)+\sum_{r=1}^{\ell-1}(\ell-r-M)\rp \right. \\
&\left. \left.  +\sum_{\ell=1}^{j-1}(j-\ell-M)          \lp \sum_{r=\ell}^{k}(r-\ell-M)+\sum_{r=1}^{\ell-1}(\ell-r-M)\rp \rb
 \left. \rule{0mm} {10mm}\rb \right/k^4
= -\frac{                     (k - 1) (k - 2) (k + 2) (k + 1) (k^2  + 5)}{3780k^3},
 \end{align*}

\bals
T_{0,1,1,1}&=\sum_{i=1}^{k}\lb \rule{0mm} {10mm}\right.\sum_{j=i}^{k}(j-i)\lb \sum_{\ell=j}^{k}(\ell-j)\lp \sum_{r=\ell}^{k}(r-\ell)
+\sum_{r=1}^{\ell-1}(\ell-r)\rp \right. \\
&\left. +\sum_{\ell=1}^{j-1}(j-\ell)\lp \sum_{r=\ell}^{k}(r-\ell)+\sum_{r=1}^{\ell-1}(\ell-r)\rp\rb\\
& \left.  +\sum_{j=1}^{i-1}(i-j)\lb \sum_{\ell=j}^{k}(\ell-j)\lp \sum_{r=\ell}^{k}(r-\ell)+\sum_{r=1}^{\ell-1}(\ell-r)\rp  
+\sum_{\ell=1}^{j-1}(j-\ell)\lp \sum_{r=\ell}^{k}(r-\ell)+\sum_{r=1}^{\ell-1}(\ell-r)\rp\rb \left. \rule{0mm} {10mm}\rb \right/k^4\\
&=\frac{                      (k - 1) (k + 1) (17 k^4  - 39 k^2  + 24)}{420k^3},
 \end{align*}

 \bals
\overline{T}_{0,1,2}&=\sum_{i=1}^{k}\lb  \sum_{j=i}^{k}(j-i-M)\lp \sum_{\ell=j}^{k}(\ell-j-M)^2
+ \sum_{\ell=1}^{j-1}(j-\ell-M)^2\rp  \right.\\
& \left. \left. +\sum_{j=1}^{i-1}(i-j-M)\lp \sum_{\ell=j}^{k}(\ell-j-M)^2+\sum_{\ell=1}^{j-1}(j-\ell-M)^2\rp \rb \right/k^3\\
&= \frac{(k - 1) (k - 2) (k + 2) (k + 1) (k^2  + 2)}{540k^3},
 \end{align*}

\bals
 T_{0,1,2}&= \sum_{i=1}^{k}\lb \sum_{j=i}^{k}(j-i)\lp \sum_{\ell=j}^{k}(\ell-j)^2+\sum_{\ell=1}^{j-1}(j-\ell)^2\rp) 
+\left. \sum_{j=1}^{i-1}(i-j)\lp \sum_{\ell=j}^{k}(\ell-j)^2+\sum_{\ell=1}^{j-1}(j-\ell)^2\rp\rb \right/k^3\\
&=\frac{ (k - 1) (k + 1) (11 k^2  - 14)}{180k},
 \end{align*}

The couple $y_i y_{i+1}$		is probabilistically reversible, hence $T_{0,1,2}=T_{0,2,1},\overline{T}_{0,1,2}=\overline{T}_{0,2,1}$.

 Now we symbolically expand (recall that $y_i$ is independent of $y_{i+2}$)
\[(\overline{y}_i+\overline{y}_{i+1}+\overline{y}_{i+2})^3,\mbox{ with }\overline{y}_i:=y_i-M.\]
Again, the contribution of $x_0,x_m$ is negligible and terms like $\overline{y}_i\overline{y}_{i+3}$ lead to $0$ by independence.
We must only retain the terms
\[S:=\overline{y}_i^3+3\overline{y}_i\overline{y}_{i+1}^2+3\overline{y}_i^2\overline{y}_{i+1}
+6\overline{y}_i\overline{y}_{i+1}\overline{y}_{i+2}.\]
Indeed, when counting the tuples, we only retain contributions  of order $m$ and neglect any other $\BO(1)$ terms or null terms (by independence).
We expand, this leads to
\bals
&(y_i^3+3 y_i  y_{i+1}^2+6  y_i  y_{i+1}  y_{i+2}+3  y_i^2  y_{i+1})+(-6  y_i^2-3  y_{i+1}^2-18  y_i  y_{i+1}-6  y_{i+1}  y_{i+2}     -6  y_i  y_{i+2}) M\\
&+(15  y_{i+1}+6  y_{i+2}+18  y_i) M^2-13 M^3.
\end{align*} 
 We make a three steps substitution, \emph{in this order}. For instance,  in $ y_i  y_{i+1}^2$, we cannot simply replace $ y_i$ by $M$ and $y_{i+1}^2$ by $T_{0,2}$. We \emph{must} use $T_{0,1,2}$.
\bit
\item $y_i^3=T_{0,3},y_i^2 y_{i+1}=T_{0,2,1}, y_i y_{i+1}^2=T_{0,1,2},y_i y_{i+1}y_{i+1}=T_{0,1,1,1},$
\item $y_i^2=T_{0,2},y_{i+1}^2=T_{0,2},y_i y_{i+1}=T_{0,1,1},  y_{i+1}y_{i+2}=T_{0,1,1},y_i y_{i+2}=M^2,$
\item $y_i=M,y_{i+1}=M,y_{i+2}=M.$
\eit
This leads to the dominant term of $\mu_3(P_n)$.

\bth
In the uniform $[1,k]$ case, the dominant term of $\mu_3(P_n)$ given by 
\bals
\mu_3(P_n)&=n\mu_3^*+\BO(1),\\
\mu_3^*&= (T_{0,3}+3T_{0,1,2}+6T_{0,1,1,1}+3T_{0,2,1})+(-9T_{0,2}-24T_{0,1,1}-6M^2)M-26M^3\\
&= \frac{4 (k - 2) (1 + 2 k) (2 k - 1) (k + 2) (k - 1) (k + 1)}{945k^3}.
\end{align*}
\ethGL  

Of course, this can also be obtained as  
\[n[\overline{ T}_{0,3}+3\overline{T}_{0,1,2}+6\overline{T}_{0,1,1,1}+3\overline{T}_{0,2,1}],\]    
but we also gave the first approach, which will be used in the next section.

The fourth centered moment $\mu_4(P_n)$ can be  similarly mechanically computed. Note that the dominant term is     there of order $n^2$: we have contribution of type $\overline{y_i}^2,\overline{y_k}^2,k\geq i+2$.

\section{The geometric$(p)$ case}\label{S4} 

We will now consider the			geometric$(p)$ case, with distribution	$pq^{i-1},i\geq 1,p\in(0,1),q:=1-p$. The computation of the 
centered cross-moments $	\overline{T_.}$ is rather intricate (in particular with many indices), even for Maple. So we will only use the ordinary 	cross-moments $	T_.$. Of course, our techniques can be applied to other polyominoes' models.

\begin{figure}[htbp]
	\centering
		\includegraphics[width=0.8\textwidth,angle=0]{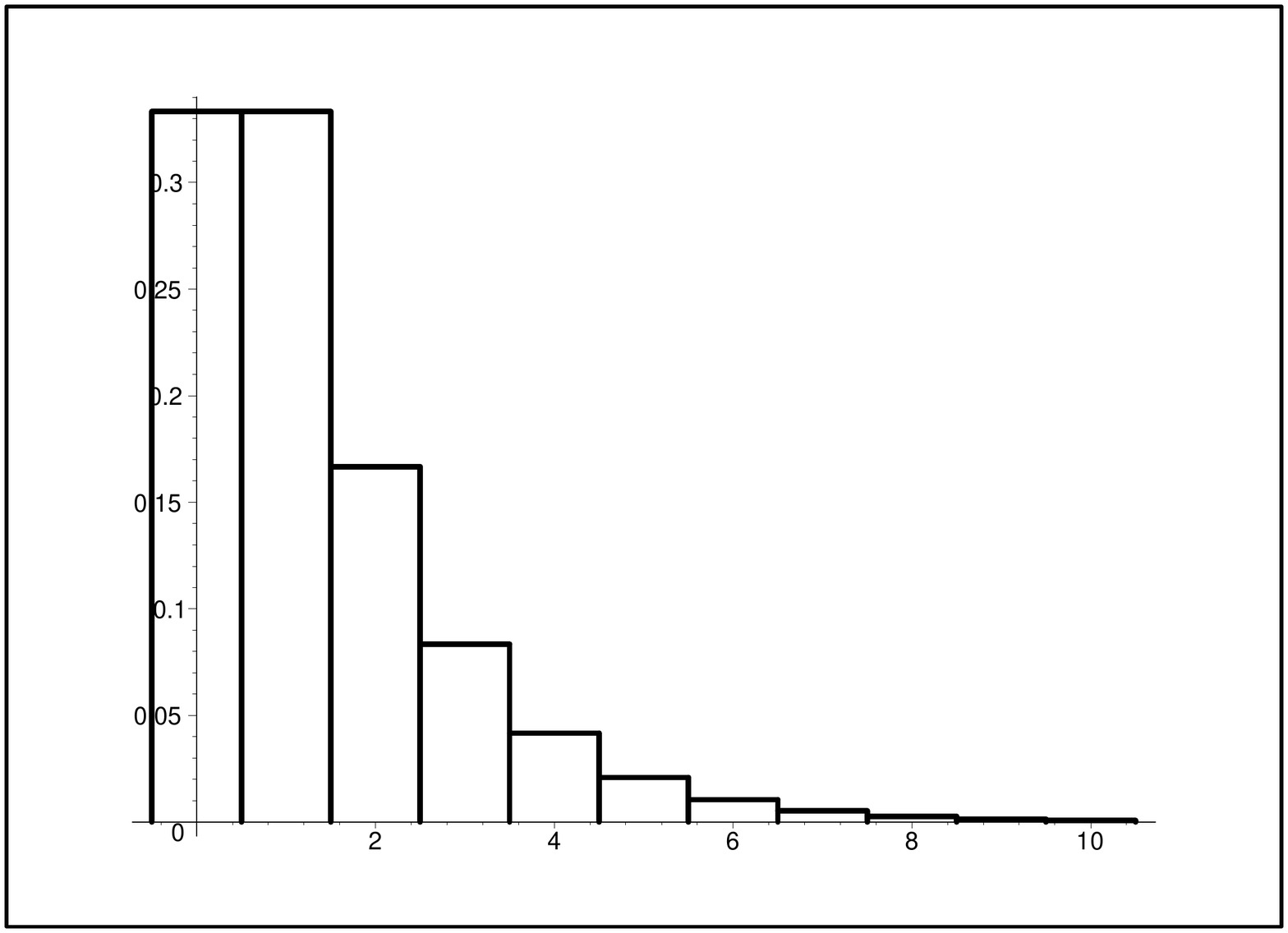}
	\caption{$f(u),p=1/2$ }
	\label{F3}
\end{figure}

 The distribution $f(u):=\P(y_i=u)$, $u$ is a non-negative integer,		is given as follows:
\bals
f(u)&=\sum_{i=1}^{\II}pq^{i-1}pq^{i+u-1}+\sum_{i=u+1}^{\II}pq^{i-1}pq^{i-u-1}=\frac{ 2 p (1 - p)^u}{2-p},u>0,\\
f(0)&=\sum_{i=1}^{\II}(pq^{i-1}pq^{i-1})=\frac{p}{2-p}.
\end{align*}  
A  plot of $f(u),p=1/2$ is given in Fig.\ref{F3}									

 The first  expressions are given as follows
\bals 
 T_{1}&=\sum_{i=1}^{\II}pq^{i-1}i=\frac{1}{p},\\
M=T_{0,1}&=\sum_{u=1}^{\II}f(u)u=\sum_{i=1}^{\II}pq^{i-1}\lp \sum_{j=i}^{\II}pq^{j-1}(j-i)+\sum_{j=1}^{i-1}pq^{j-1}(i-j)\rp
=\frac{ 2 (1-p)}{ p (2-p)},\\
 \end{align*}  
hence
\bals
 M_{R,n}&=(n-1)M+2T_1=\frac{ 2 + (2 - 2 p) n}{ p (2-p)},\\ 
M_{P,n}&=(n-1)M+2n+2T_1=\frac{ 2 + (2 + 2 p - 2 p^2 ) n}{p (2-p)}.
 \end{align*}

We recall a previous definition:
\bals
T_{\al,\be, \gam,\de}&:=\E\lp x_0^\al\cdot y_1^\be \cdot y_2^\gam \cdot y_3^\de \rp\\
 &=\left. \lb \sum_{i=1}^{k}i^\al \sum_{j=1}^{k}|j-i|^\be \sum_{\ell=1}^{k}|\ell-j|^\gam \sum_{r=1}^{k}|r-\ell|^\de\rb\right/k^4 .
 \end{align*} 
													The next necessary expressions are given as follows:	
\bals
 T_{0,2}&=\sum_{i=1}^{\II}pq^{i-1}\lp \sum_{j=i}^{\II}pq^{j-1}(j-i)^2 +\sum_{j=1}^{i-1}pq^{j-1}(i-j)^2\rp 
=\sum_{u=1}^{\II}f(u)u^2= \frac{     2 (1-p)}{p^2},\\
 T_{0,3}&=\sum_{u=1}^{\II}f(u)u^3=\sum_{i=1}^{\II}pq^{i-1}\lp \sum_{j=i}^{\II}pq^{j-1}(j-i)^3
 +\sum_{j=1}^{i-1}pq^{j-1}(i-j)^3\rp \\
&=\frac{2 (1-p) (p^2  - 6 p + 6)}{     p^3  (2-p)},\\
 T_{0,1,1}&=\sum_{i=1}^{\II}pq^{i-1}\lb \sum_{j=i}^{\II}pq^{j-1}(j-i)\lp \sum_{\ell=j}^{\II}pq^{\ell-1}(\ell-j) 
+\sum_{\ell=1}^{j-1}pq^{\ell-1}(j-\ell)\rp \right.\\
&\left. +\sum_{j=1}^{i-1}pq^{j-1}(i-j)\lp \sum_{\ell=j}^{\II}pq^{\ell-1}(\ell-j)+\sum_{\ell=1}^{j-1}pq^{\ell-1}(j-\ell)\rp\rb\\
&= \frac{                   (1-p) (p^4  - 7 p^3  + 23 p^2  - 32 p + 16)}{     p^2  (2-p)^2  (p^2  + 3 - 3 p)},
\end{align*}

\bals
 T_{0,1,1,1}&=\sum_{i=1}^{\II}pq^{i-1}\lb \rule{0mm} {10mm}\right.\sum_{j=i}^{\II}pq^{j-1}(j-i)\lb \sum_{\ell=j}^{\II}pq^{\ell-1}(\ell-j) 
\lp \sum_{r=\ell}^{\II}pq^{r-1}(r-\ell) +\sum_{r=1}^{\ell-1}pq^{r-1}(\ell-r)\rp \right. \\
& \left. +\sum_{\ell=1}^{j-1}pq^{\ell-1}(j-\ell) \lp \sum_{r=\ell}^{\II}pq^{r-1}(r-\ell)+\sum_{r=1}^{\ell-1}pq^{r-1}(\ell-r)\rp\rb\\
& +  \sum_{j=1}^{i-1}pq^{j-1}(i-j)\lb \sum_{\ell=j}^{\II}pq^{\ell-1}(\ell-j)     \lp \sum_{r=\ell}^{\II}pq^{r-1}(r-\ell)\right. \right.   \\                               &\left. \left.  +\sum_{r=1}^{\ell-1}pq^{r-1}(\ell-r)\rp   +\sum_{\ell=1}^{j-1}pq^{\ell-1}(j-\ell)\lp \sum_{r=\ell}^{\II}pq^{r-1}(r-\ell)+\sum_{r=1}^{\ell-1}pq^{r-1}(\ell-r)      \rp  \rb\left.   \rule{0mm} {10mm}\rb\\
&= \frac{2 (28-84p +113 p^2  - 86 p ^3 + 39 p ^4 - 10 p^5  + p^6 ) (1-p) ^2}{p^3  (p^2  - 2 p + 2) (2-p) (p^2  + 3 - 3 p)^2 }, 
\end{align*}

\bals
 T_{0,1,2}&=\sum_{i=1}^{\II}pq^{i-1}\lb  \sum_{j=i}^{\II}pq^{j-1}(j-i)\lp \sum_{\ell=j}^{\II}pq^{\ell-1}(\ell-j)^2             
           +\sum_{\ell=1}^{j-1}pq^{\ell-1}(j-\ell)^2\rp\right.\\
					&\left. +\sum_{j=1}^{i-1}pq^{j-1}(i-j) \lp \sum_{\ell=j}^{\II}pq^{\ell-1}(\ell-j)^2
					+\sum_{\ell=1}^{j-1}pq^{\ell-1}(j-\ell)^2\rp\rb      \\
&=\frac{                  (28 - 56 p + 38 p^2  - 10 p^3  + p^4 ) (1-p)}{    p^3  (2-p)^3}, 
\end{align*}

 \[T_{1,1}=\sum_{i=1}^{\II}pq^{i-1}i\lp \sum_{j=i}^{\II}pq^{j-1}(j-i)+\sum_{j=1}^{i-1}pq^{j-1}(i-j)\rp 
= \frac{    (1-p) (p^2  - 4 p + 6)}{   p^2  (2-p)^2}.\]

 Again, $T_{0,1,2}=T_{0,2,1}$.

	The dominant term of $\V(R_n)$ is given by \\		(all necessary expressions are extracted from Sec. \ref{S2} and \ref{S3})												
  \bals                             
  &n[(T_{0,2}-M^2)+2(T_{0,1,1}-M^2)]=nV^*,\\
	V^*&=\frac {4 (1-p) (p^4  + 9 p^2  - 4 p^3  - 10 p + 5)}{p ^2 (2-p)^2  (p^2  + 3 - 3 p)} . 
	\end{align*}

 The exact value of $\V(P_n)$ is given by                               

 \bals
 \V(R_n)&=(n-1)T_{0,2}+2T_{2}+(n-2)2T_{0,1,1}+4T_{1,1}+(n-2)(n-3)M^2+2T_1^2+4T_1M(n-2)-M_{R,n}^2\\
&=\frac{n[ 4 (1-p) (p ^4 + 9 p^2  - 4 p ^3 - 10 p + 5)]+ 4 (3 p ^2 - 5 p + 5) (1-p)^2}{p^2  (2-p)^2 (p^2  + 3 - 3 p)}.
\end{align*}

The third centered moment $\mu_3(P_n)$ (dominant term)  is given by 

  \bals
\mu_3(P_n)&=n[ (T_{0,3}+3T_{0,1,2}+6T_{0,1,1,1}+3T_{0,2,1})+(-9T_{0,2}-24T_{0,1,1}-6M^2)M-26M^3]+\BO(1)\\
&= n\frac{  8 (1-p) (114 - 570 p + 1332 p^2  - 1908 p^3  + 1849 p^4            
         - 1263 p^5  + 616 p^6  - 213 p^7  + 52 p^8  - 9 p^9 + p^{10})}{(2-p)^3 p^3 (p^2  - 2 p + 2) (p^2  + 3 - 3 p)^2}.\\
				&+\BO(1).
\end{align*}

We summarize our results in the following theorem
\bth
The first three moments of $P_n$ in the geometric$(p)$ case are given by
\bals
M_{P,n}&=(n-1)M+2n+2T_1=\frac{- 2 + (-2 - 2 p + 2 p^2 ) n}{p (-2 + p)},\\
\V(P_n)&=\frac{n[ 4 (1-p) (p ^4 + 9 p^2  - 4 p ^3 - 10 p + 5)]+ 4 (3 p ^2 - 5 p + 5) (1-p)^2}{p^2  (2-p)^2 (p^2  + 3 - 3 p)},\\
\mu_3(P_n)&= n\frac{  8 (1-p) (114 - 570 p + 1332 p^2  - 1908 p^3  + 1849 p^4            
         - 1263 p^5  + 616 p^6  - 213 p^7  + 52 p^8  - 9 p^9 + p^{10})}{(2-p)^3 p^3 (p^2  - 2 p + 2) (p^2  + 3 - 3 p)^2}\\
				&+\BO(1).
\end{align*}  
\ethGL

\section{The stochastic processes in the uniform $[1,k]$ case}\label{S5} 
In this section, we analyze  the stochastic processes related to $P_n$. Seen as a stochastic process, the random part of the perimeter is asymptotically given
by $Q_m(j):=\sum_i^j y_i$: we can ignore $x_0,x_m$ and the contribution $2n$ is a constant. By the functional central limit theorem 
  (\cite[p.~174, Thm. 20.1]{BI68}),
we obtain the following result, where $B(t)$ is the standard Brownian 
Motion (BM) and $\Rightarrow $ denotes the weak convergence of random functions in the space of all right-continuous functions 
that have right limits and are endowed with the Skorohod metric (the $\Fi$-mixing property is immediate here: see (\cite[p.~167,     example $1$]{BI68})). This gives the limiting trajectories corresponding to $Q_m(j)$.
\bth
\[\frac{Q_m (\lfloor mt\rfloor ) - M  mt}{\sig \sqrt{m}}
\Rightarrow B(t),\B
m\rightarrow \infty,\  t \in [0,1],\sig=\sqrt{V^*}.\]
\ethGL
As a corollary, we have
\bth                                                           \label{Thm1}
\[\frac{Q_m-mM}{\sig \sqrt{m}}\sim \BN(0,1),m\rightarrow \infty,\]                       
\ethGL
where $\BN$ is a  Gaussian (normal) random variable.

\begin{figure}[htbp]
	\centering
		\includegraphics[width=0.8\textwidth,angle=0]{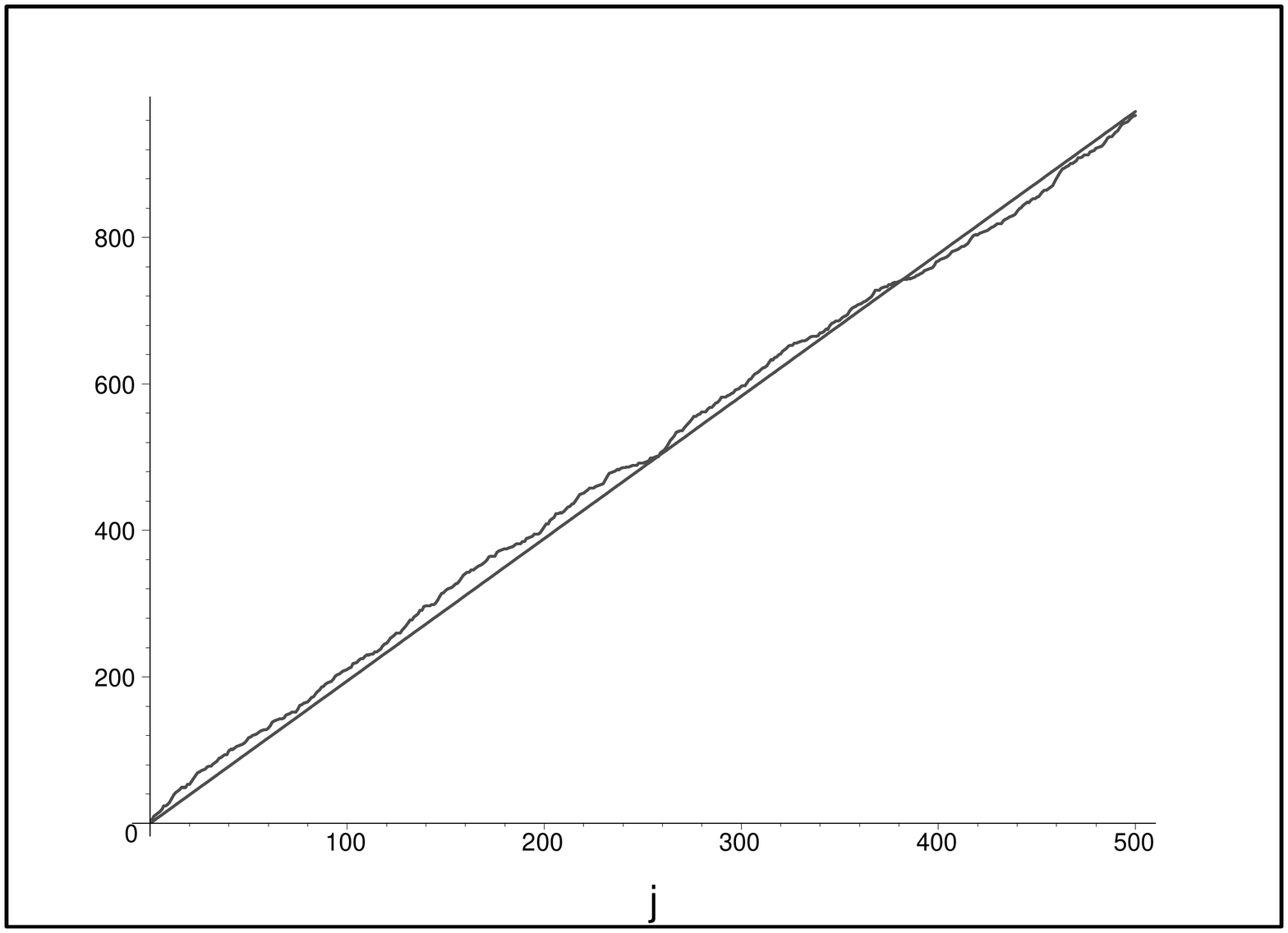}
	\caption{$Q_m(j),m=500$, drift= $jM$ }
	\label{F4}
\end{figure}

\begin{figure}[htbp]
	\centering
		\includegraphics[width=0.8\textwidth,angle=0]{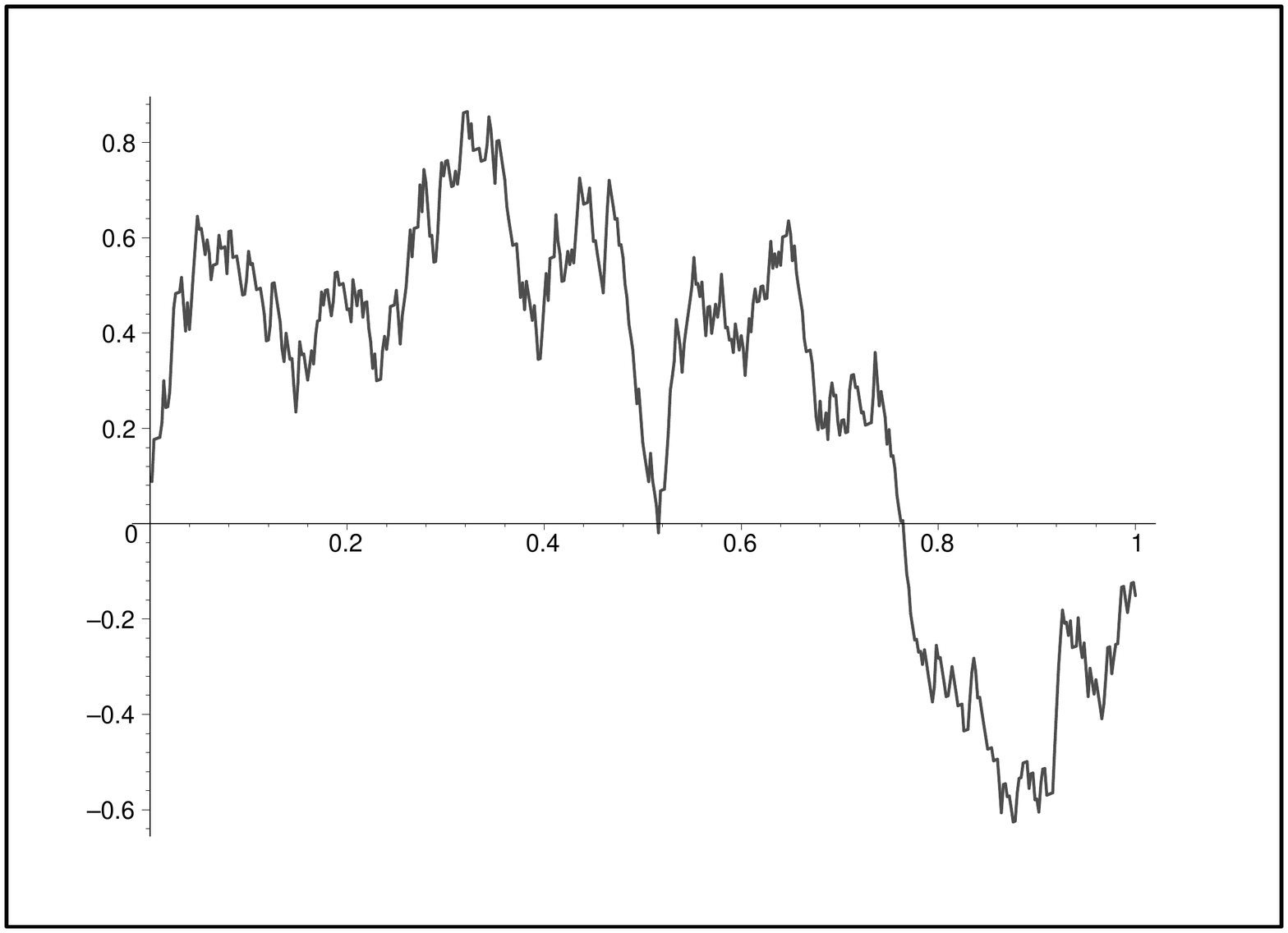}
	\caption{a typical normalized trajectory }
	\label{F5}
\end{figure}

In the uniform case, $k=6$, we have made a simulation of $N=100000$ trajectories $Q_m(j),m=500$. A typical trajectory is given in Fig.\ref{F4}, together with the drift $jM$.
In Fig.\ref{F5}, we show a typical normalized trajectory	
\[\frac{Q_m (\lfloor mt\rfloor ) - M  mt}{\sig \sqrt{m}}	,\]
with the classical strongly irregular $BM$ behaviour. 
We have also computed the observed moments: let $z_\ell$ denote the $ \ell$th simulated value of $Q_m(m)-mM$. We obtain
\[\left.  \lp\sum_{\ell=1}^N  \frac{z_\ell}{\sig \sqrt{m}}\rp\right/ N=-0.0038\ldots,\left. \lp\sum_{\ell=1}^N  \lb\frac{z_\ell}{\sig \sqrt{m}}\rb^2\rp \right/N=0.991\ldots,
\sum_{\ell=1}^N  z_\ell^3= 1287.47,\]		
to be compared with the theoretical values $\{0,1,m\mu_3^*=1569.272976\ldots\}$.
About the third moment, another simulation gives $1911.44\ldots$: $m$ is not large enough to give a really good fit.
 
\begin{figure}[htbp]
	\centering
		\includegraphics[width=0.8\textwidth,angle=0]{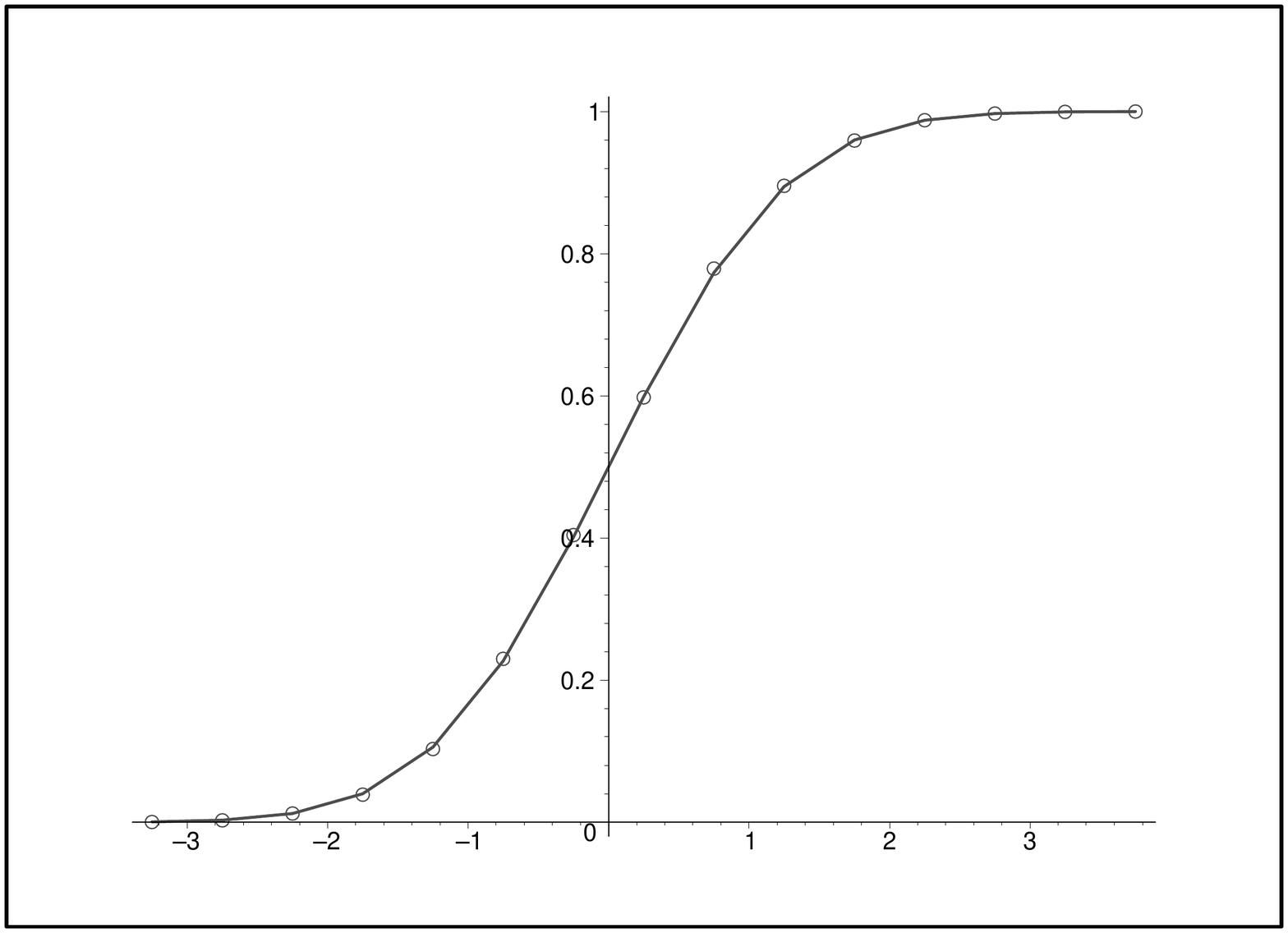}
	\caption{the cumulative histogram (circle)  and the Gaussian distribution function (line) }
	\label{F13}
\end{figure}

To illustrate Thm \ref{Thm1}, we have build a histogram as follows: we construct a set of intervals
 $I(i):=[i\De-3-3\De/2,i\De-3-\De/2],i=0..(6/\De+2)$, centered on $i\De-3-\De$ and covering the interval $[-3-\De,3+\De]$. We choose here $\De=1/2$. We define cells such that $cell(i)$ corresponds to interval $I(i)$. We compute the number $N[i]$ of values of $\frac{z_\ell}{\sig \sqrt{m}}$
falling into interval $I(i)$ and put $N(i)/N$ into $cell(i)$ (values $<3.5$ are attributed to $cell(0)$ and similarly for values $>3.5$). This gives the empirical histogram. In Fig.\ref{F13}, we compare the cumulative histogram (circle)  with the Gaussian distribution function( line): the fit is quite good. 

\begin{figure}[htbp]
	\centering
		\includegraphics[width=0.8\textwidth,angle=0]{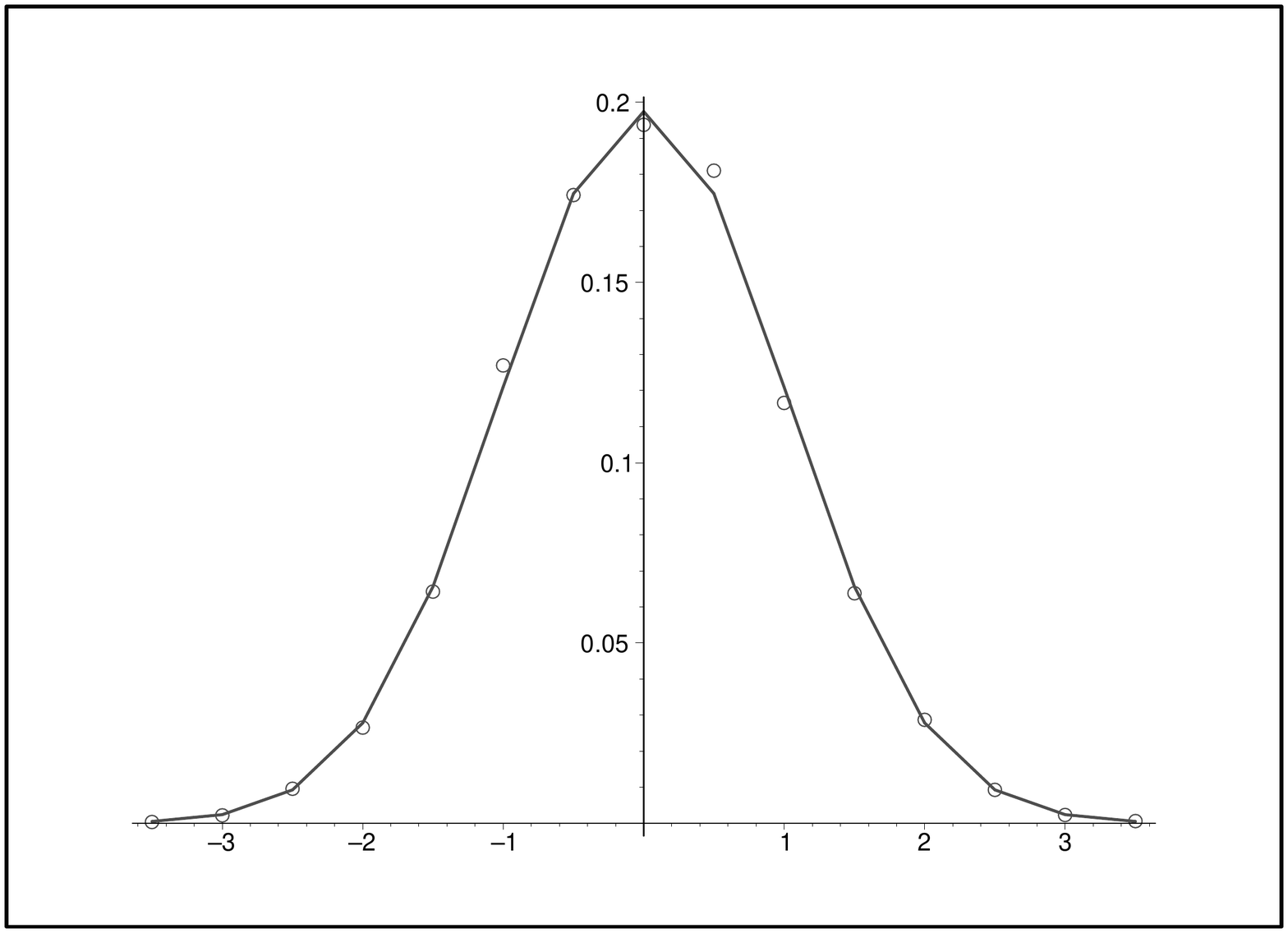}
	\caption{the histogram (circle) and the Gaussian probability mass function in each interval $I(i)$ (line) }
	\label{F11}
\end{figure}

But it it still more precise to compare, in Fig.\ref{F11} the histogram itself (circle) with the Gaussian probability mass in interval $I(i)$: 
$\int_{i\De-3-3\De/2}^{i\De-3-\De/2}\exp(-x^2/2)/\sqrt{2\pi}dx$ (line). The fit is quite satisfactory.

We have also made the same kind of simulations for the geometric$(p)$  case. The results are quite similar. 

\section{Conclusion}\label{S6} 
We have shown that a probabilistic approach leads, almost mechanically, to  the first three moments of $P_n$ and its asymptotic Brownian and Gaussian properties.  This technique can be applied to other moments and to other initial probability distributions.

\bibliographystyle{plain} 

\end{document}